\documentclass[12pt,a4paper,reqno]{amsart}

\usepackage{graphicx}
\usepackage{amsfonts,amsthm,latexsym,amsmath,amssymb,amscd,epsfig}
\usepackage{graphics,graphicx}
\usepackage[all]{xy}
\usepackage{young}

\usepackage{tikz}
\usetikzlibrary{matrix}

\newtheorem{theorem}{Theorem}[section]
\newtheorem{lemma}[theorem]{Lemma}

\theoremstyle{definition}

\newtheorem{proposition}[theorem]{Proposition}
\newtheorem{corollary}[theorem]{Corollary}
\theoremstyle{remark}

\numberwithin{equation}{section}

\newcommand{\nc}{\newcommand}
\renewcommand{\frak}{\mathfrak}
\providecommand{\cal}{\mathcal}
\renewcommand{\bold}{\mathbf}

\newcommand{\ZZ}{{\mathbb Z}}

\nc \Ab{{\ensuremath{\bold A}}}
\nc \ab{{\ensuremath{\bold a}}}
\nc \bb{{\ensuremath{\bold b}}}
\nc \cb{{\ensuremath{\bold c}}}
\nc \Bb{{\ensuremath{\bold B}}}
\nc \Gb{{\ensuremath{\bold G}}}
\nc \Qb{{\ensuremath{\bold Q}}}
\nc \Rb{{\ensuremath{\bold R}}} \nc \Cb{{\ensuremath{\bold C}}} 
\nc \Eb{{\ensuremath{\bold E}}}
\nc \eb{{\ensuremath{\bold e}}}
\nc \Db{{\ensuremath{\bold D}}}
\nc \Fb{{\ensuremath{\bold F}}}
\nc \ib{{\ensuremath{\bold i}}}
\nc \jb{{\ensuremath{\bold j}}}
\nc \kb{{\ensuremath{\bold k}}}
\nc \nb{{\ensuremath{\bold n}}}
\nc \rb{{\ensuremath{\bold r}}}
\nc \Pb{{\ensuremath{\bold P}}}
\nc \pb{{\ensuremath{\bold p}}}
\nc \SPb{{\ensuremath{\bold {SP}}}}
\nc \Zb{{\ensuremath{\bold Z}}} 
\nc \zb{{\ensuremath{\bold z}}} 
\nc \gb{{\ensuremath{\bold g}}} 
\nc \fb{{\ensuremath{\bold f}}} 
\nc \ub{{\ensuremath{\bold u}}} 
\nc \vb{{\ensuremath{\bold v}}} 
\nc \yb{{\ensuremath{\bold y}}} 
\nc \xb{{\ensuremath{\bold x}}} 
\nc \xib{{\ensuremath{\bold \xi}}} 
\nc \Nb{{\ensuremath{\bold N}}} 
\nc \Hb{{\ensuremath{\bold H}}} 
\nc \wb{{\ensuremath{\bold w}}} 
\nc \Wb{{\ensuremath{\bold W}}} 
\nc \syz{{\mathbf {syz}}}
\nc \bnoll{{\ensuremath{\bold 0}}} 

\nc \mf{\frak m} \nc \mh{\hat{\m}} 
\nc \nf{\frak n}
\nc \Of{\frak O}
\nc \rf{\frak r}
\nc \mufr{{\mathbf \mu}}
\nc \hf{\frak h} 
\nc \qf{\frak q} 
\nc \bfr{\frak b} 
\nc \kfr{\frak k} 
\nc \pfr{\frak p} 
\nc \af{\frak a }
\nc \cf{\frak c }
\nc \sfr{\frak s} 
\nc \ufr{\frak u} 
\nc \g{\frak g} 
\nc \gA{\g_{\Ao}} 
\nc \lfr{\frak l}
\nc \afr{\frak a}
\nc \gfh{\hat {\frak g}}
\nc \gl{\frak { gl }}
\nc \Sl{\frak {sl}}
\nc \SU{\frak {SU}}
\nc{\Homf}{\frak{Hom}}

\newcommand{\on}{\operatorname}
\nc\hankel{\on {Hankel}}
\nc\row{\on {row\ }}
\nc\nullity{\on {nullity }}
\nc\col{\on {col\ }}
\nc\rowm{\on {Row \ }}
\nc\loc{\on {lc \ }}
\nc\nullo{\on {null\ }}
\nc\Nul{\on {Nul\ }}
\nc \Ann {\on {Ann }}
\nc \Ass {\on {Ass \ }}
\nc \Coker {\on {Coker}}
\nc \Co{\on C}
\nc \Homo{\on {Hom}}
\nc \Ker {\on {Ker}}
\nc \omod{\on {mod}}
\nc \No {\on N}
\nc \NN {\on {NN}}
\nc \NGo {\on {NG}}
\nc \Oo {\on O}
\nc \ch {\on {ch}}
\nc \rko {\on {rk}}
\nc \Sing {\on {Sing\ }}
\nc \Reg {\on {Reg}}
\nc \CoI {\on {CI}}
\nc \CoM {\on {CM}}
\nc \Gor {\on {Gor}}
\nc \Type {\on {Type}}
\nc \can {\on {can}}
\nc \Top {\on {T}}
\nc \rel {\on {rel}}
\nc \sgn {\on {sgn }}
\nc \trdeg {\on {tr.deg}}
\nc \codim {\on {codim }}
\nc \coht {\on {coht}}
\nc \divo {\on {div \ }}
\nc \coh {\on {coh}}
\nc \Clo {\on {Cl}}
\nc \embdim{\on {embdim}}
\nc \embcodim{\on {embcodim \ }}
\nc \qcoh {\on {qcoh}}
\nc \grad {\on {grad}\ }
\nc \grade {\on {grade}}
\nc \hto {\on {ht}}
\nc \depth {\on {depth}}
\nc \prof {\on {prof}}
\nc \reso{\on {res}}
\nc \ind{\on {ind}}
\nc \prodo{\on {prod}}
\nc \coind{\on {coind}}
\nc \Con{\on {Con}}
\nc \Crit{\on {Crit}}
\nc \Der{\on {Der}}
\nc \Char{\on {Char}}
\nc \Ch{\on {Ch}}

\nc \Ext{\on {Ext}}
\nc \Eo{\on {E}}
\nc \End{\on {End}}
\nc \Ad{\on {Ad}}
\nc \gr{\on {gr}}
\nc \Fo{\on {F}}
\nc \Gr{\on {Gr}}
\nc \Go{\on {G}}
\nc \GFo{\on {GF}}
\nc \Glo{\on {Gl}}
\nc \Ho{\on {H}}
\nc \CMo{\on {\CM}}
\nc \SCM{\on {SCM}}
\nc \hol{\on {hol}}
\nc{\sgd}{\on{sgd}}
\nc \supp{\on {supp}}
\nc \ssupp{\on {s-supp}}
\nc \singsupp{\on {singsupp}}
\nc \msupp{\on {msupp}}
\nc \spec{\on {spec}}
\nc \spano{\on {span }}
\nc \Span{\on {Span }}
\nc \Max{\on {Max}}
\nc \Min{\on {Min}}
\nc \Mod{\on {Mod}}
\nc \Rad {\on {Rad}}
\nc \rad {\on {rad}}
\nc \rank {\on {rank}}
\nc \range {\on {range}}
\nc \Slo{\on {SL}}
\nc \soc {\on {soc}}
\nc \Irr {\on {Irr}}
\nc \Imo {\on {Im}}
\nc \SSo{\on {SS}}
\nc \lub{\on {lub}}
\nc \gldim{\on {gl.d.}}
\nc \pdo{\on {p.d.}} 
\nc \ido{\on {i.d.}} 
\nc \dSSo{\dot {\SSo}}
\nc \So{\on S}
\nc \Io{\on I}
\nc \Jo{\on J}
\nc \jo{\on j}
\nc \Ko{\on K}
\nc \PBW{\Ac_{PBW}}
\nc \Ro{\on R}
\nc \To{\on T}
\nc \Ao{\on A}

\nc \Do{{\on D}}
\nc \Bo{\on B}
\nc \Po{\on P}
\nc \Qo{\on Q}
\nc \Zo{\on Z}
\nc \U{\on U}
\nc \wt{\on {wt}}
\nc \Uh{\hat {\U}}
\nc \T{\on T}
\nc \Lo{\on L}
\nc{\dop}{\on d}
\nc{\eo}{\on e}
\nc{\ado}{\on{ad}}
\nc{\Tot}{\on{Tot}}
\nc{\Aut}{\on{Aut}}
\nc{\sinc}{\on {sinc}}

\nc{\overrightleftarrows}[2]{\overset{#1}{\underset{#2}{\rightleftarrows}}}

\nc{\CCF}{\cal{CF}}

\nc{\CDF}{\cal{DF}}
\nc{\CHC}{\check{\cal C}}

\nc{\Cone}{\on{Cone}}
\nc{\dec}{\on{dec}}

\nc{\Diff}{\on{Diff}}

\nc{\dirlim}{\underset{\to}{\on{lim}}}
\nc{\dpar}{\partial}
\nc{\GL}{\on{GL}}
\nc{\CGr}{\cal{G}r}

\nc{\pr}{\on{pr}}
\nc{\semid}{|\!\!\!\times}
\nc{\Hom}{\on{Hom}}
\nc \RHom{\on {RHom}}

\nc \Proj{\mathrm {Proj\ }}
\nc \proj{\mathrm {proj}}
\nc{\Id}{\on{Id}}
\nc{\id}{\on{id}}
\nc{\Ima}{\on{Im}}
\nc{\invtimes}{\underset{\gets}{\otimes}}
\nc{\invlim}{\underset{\gets}{\on{lim}}}
\nc{\Lie}{\on{Lie}}

\nc{\re}{\on{Re }}
\nc{\Pic}{\on{Pic }}
\nc{\LPic}{\on{LPic }}
\nc{\Sch}{\on{Sch}}
\nc{\Sh}{\on{Sh}}
\nc{\Set}{\on{Set}}
\nc{\spo}{\on{sp\  }}
\nc{\Spec}{\on{Spec}}
\nc{\mSpec}{\on{mSpec}}
\nc{\Specb}{\bold {Spec}}
\nc{\Projb}{\bold {Proj}}
\nc{\Specan}{\on{Specan}}
\nc{\Spo}{\on{Sp}}
\nc{\Spf}{\on{Spf}}
\nc{\sym}{\on{sym}}
\nc{\symm}{\on{symm}}
\nc{\rop}{\on{r}}
\nc{\Td}{\on{Td}}
\nc{\Tor}{\on{Tor}}

\nc{\Artin}{\cal{A}rtin}
\nc{\Dgcoalg}{\cal{D}gcoalg}
\nc{\Dglie}{\cal{D}glie}
\nc{\Ens}{\cal{E}ns}
\nc{\Fsch}{\cal{F}sch}
\nc{\Groupoids}{\cal{G}roupoids}
\nc{\Holie}{\cal{H}olie}
\nc{\Mor}{\cal{M}or}

\nc{\CF}{\ensuremath{\cal{F}}}
\nc \Kc{\ensuremath{\cal K}}
\nc \Lc{{\ensuremath{\cal L}}}
\nc \lcc{{\mathcal l}} 
\nc \CC{{\ensuremath{\cal C}}} 
\nc \Cc{{\ensuremath {\cal C}}}
\nc \Pc{{\ensuremath{\cal P}}}
\nc \Dc{\ensuremath{\mathcal D}}
\nc \Ac{{\ensuremath{\cal A}}} 
\nc \Bc{{\ensuremath{\cal B}}}
\nc \Ec{{\ensuremath{\cal E}}}
\nc \Fc{{\ensuremath{\cal F}}}
\nc \Mcc{{\ensuremath{\cal M}}} 
\nc \hM{\hat{\Mcc}} 
\nc \bM{\bar {\Mcc}} 
\nc\hbM{\hat{\bar \Mcc}}  
\nc \Nc{{\ensuremath{\cal N}}}
\nc \Hc{{\ensuremath{\cal H}}} 
\nc \Ic{{\ensuremath{\cal I}}} 
\nc \Oc{\ensuremath{{\cal O}}}
\nc \Och{\hat{\cal O}} 
\nc \Sc{{\ensuremath{{\cal S}}}}
\nc \Tc{\ensuremath{{\cal T}}} 
\nc \Vc{{\ensuremath{{\cal V}}}} 
\nc{\CA}{{\ensuremath{{\cal A}}}}
\nc{\CB}{{\ensuremath{{\cal B}}}}
\nc{\C}{{\ensuremath{{\cal F}}}}
\nc{\Gc}{{\ensuremath{{\cal G}}}}
\nc{\CH}{\ensuremath{\mathcal H}}
\nc{\CI}{{\ensuremath{{\cal I}}}}
\nc{\CM}{{\ensuremath{{\cal M}}}}
\nc{\CN}{{\ensuremath{{\cal N}}}}
\nc{\CO}{{\ensuremath{{\cal O}}}}
\nc{\Rc}{{\ensuremath{{\cal R}}}}
\nc{\CT}{{\ensuremath{\mathcal T}}}
\nc{\CU}{\ensuremath{{\cal U}}}
\nc{\CV}{\ensuremath{{\cal V}}}
\nc{\CZ}{\ensuremath{{\cal Z}}}
\nc{\Homc}{\ensuremath{{\cal {Hom}}}}

\nc{\Tab}{\ensuremath{{\mbox{Tab}}}}

\nc{\STab}{\ensuremath{{\mbox{STab}}}}

\nc{\fa}{\frak{a}}
\nc{\fA}{\frak{A}}
\nc{\fg}{\frak{g}}
\nc{\fh}{\frak{h}}
\nc{\fI}{\frak{I}}
\nc{\fK}{\frak{K}}
\nc{\fm}{\frak{m}}
\nc{\fP}{\frak{P}}
\nc{\fS}{\frak{S}}
\nc{\ft}{\frak{t}}
\nc{\fX}{\frak{X}}
\nc{\fY}{\frak{Y}}

\nc{\bF}{\bar{F}}
\nc{\bCP}{\bar{\cal{P}}}
\nc{\bm}{\mbox{\bf{m}}}
\nc{\bT}{\mbox{\bf{T}}}
\nc{\bS}{\mbox{\bf{S}}}
\nc{\hB}{\hat{B}}
\nc{\hC}{\hat{C}}
\nc{\hP}{\hat{P}}
\nc{\htest}{\hat P}

\nc{\nen}{\newenvironment}
\nc{\ol}{\overline}
\nc{\ul}{\underline}
\nc{\ra}{\to}
\nc{\lla}{\longleftarrow}
\nc{\lra}{\longrightarrow}
\nc{\Lra}{\Longrightarrow}
\nc{\Lla}{\Longleftarrow}
\nc{\Llra}{\Longleftrightarrow}
\nc{\hra}{\hookrightarrow}
\nc{\iso}{\overset{\sim}{\lra}}

\nc{\dsize}{\displaystyle}
\nc{\sst}{\scriptstyle}
\nc{\tsize}{\textstyle}

\nen{exa}[1]{\label{#1}{\bf Example.\ } }{}

\nen{rem}[1]{\label{#1}{\em Remark.\ } }{}
\nen{exer}[1]{\label{#1}{\em Exercise.\ } }{}

\newcommand{\la}{\lambda}

\newcommand {\bC} {\mathbb C}
\newcommand {\bR} {\mathbb R}

\newcommand{\bN}{\mathbb N}
\newcommand {\D} {\mathcal D}

\newcommand {\eps} {\epsilon}

\newcommand{\bla}{{\bf\lambda}}
\newcommand{\DD}{\tilde{\mathcal{D}}_Y}

\newcommand{\OO}{\tilde{\mathcal{O}}_X}

\begin{document}
\title{Differential operators and reflection group of type $B_n$}

\author{Ibrahim Nonkan\'e, Lat\'evi M. Lawson}

\address{Departement d'\'economie et de math\'ematiques appliqu\'ees, IUFIC, Universit\'e Thomas Sankara, Burkina faso}
\email{inonkane@univ-ouaga2.bf}

\address{ African Institute for Mathematical Sciences (AIMS), Summerhill Estates, East Legon Hills, Santoe, Accra,  Ghana P.O. Box LG DTD 20046, Legon, Accra, Ghana}
\address{  D\'epartement de math\'ematiques, Facult\'es des Sciences, Universit\'e de Lom\'e, 01 BP 1515, Lom\'e, Togo.}
\email{latevi@aims.edu.gh}

\maketitle
\begin{abstract}
In this note, we study the polynomial representation of the quantum Olshanetsky-Perelomov system for a finite reflection group $W$ of type $B_n$. We endow the polynomial ring  ${\mathbb C} [x_1,\ldots\\\ldots, x_n]$ with a structure  of module over the Weyl algebra associated with the ring  ${\mathbb C} [x_1,\ldots,x_n]^{W}$ of invariant polynomials under a reflections group $W$  of type $B_n$. Then we study the polynomial representation of the ring of invariant differential operators under the reflections group $W$. We use  the group representation theory namely the higher Specht polynomials associated with the reflection group $W$ and establish a decomposition  of  that structure by providing explicitly the generators of the simple components.\\
 
\textbf{keywords}:\ { Differential operators, Representation theory, Higher Specht polynomials,  Olshanetsky-Perelomov system,  Primitive idempotents, Reflection groups, Young diagram.}
\\\\
\textbf{ Mathematics Subject Classification}:\ { Primary 13N10, Secondary 20C30}.
\end{abstract}

\section{Introduction} 
The  Weyl algebra's history has begun with  the birth of quantum mechanics. Group theory has played a major role in the discovery of  general laws of quantum theory. It is not surprising that concepts arising in  group theory find their applications in physics.  As Hermann  Weyl has said in \cite{Weyl}, there exists  a plainly discernible parallelism between the more recent developments of mathematics and physics. Also, the theory of group representations is one of the best example of interaction between fundamental physics and pure mathematics. Having this at hand we  study  the polynomials representation of the Weyl algebra with special consideration for physical problems.  
Since  symmetries are relevant in  physics, our attention has been drawn particularly to the reflection  group $W(B_n)$  of type $B_n$ and its representation in connection with rational Olshanetsky-Perelomov operators. In what follows we study the action of invariant differential operators  on the polynomial ring through the representation theory  of reflection  groups  of type $B_n$. We show that the polynomial representation of the rational Olshanetsky-Perelomov systems for $W$  is related to the  representation theory  of $W$. In fact we have studied the polynomials representation of the rational quantum Calogero-Moser system in \cite{DOS} and similar topic in \cite{RTG, IBK}. Since the rational  Olshanetsky-Perelomov system is a generalization of rational quantum Calogero-Moser system, we extend the obtained results  in \cite{DOS} to the rational Olshanetsky-Perelomov system for reflection group of type $B_n$.

\section{Preliminaries and motivations}

\subsection{Reflection groups and root systems}
We recall some basic facts about real reflection groups of type $B_n$ (see \cite{Hu} for more details).\\
Let $\fh$ be a (real) euclidean space endowed with a positive definite symmetric linear form  $(\cdot, \cdot ).$
A reflection is a bilinear operator $s$ on $\fh$ which sends some nonzero vector $\alpha$ to its negative while fixing pointwise the hyperplane $H_{\alpha}$ orthogonal to $\alpha$. We may write $s=s_{\alpha}$ and there a simple formula: 
$$ s_\alpha = \lambda - \frac{2 (\lambda, \alpha)}{(\alpha, \alpha)},\ \forall \lambda \in \fh.$$
 The operator $s_\alpha$ is element of order $2$ in the groupe $O(\hf)$ of all orthogonal transformations of $\fh$.
 
 A finite  group generated by reflections is called a  {\it reflection group}, such group is a subgroup  of $O(\fh)$.
 Let  symmetric group $\Sc_n$ acts on $\fh=\bR^n$ by permuting the standard basis $\eps_1, \ldots,\eps_n$ of $\bR^n$. The transposition$(i,j)$ acts as a reflection, sending $\eps_i-\eps_j$ to its negative and fixing pointwise the orthogonal complement, which consists of all vectors in $\fh$ having equal $i$th and $j$th components. It is clear that $\Sc_n$ is  a reflection group. Other reflections can be made by sending an $\eps_i$ to its  negative and fixing all other $\eps_j.$ These sign changes generate a group of order $2^n$ isomorphic to $(\ZZ / 2 \ZZ)^n$, which intersects $\Sc_n$ trivially and is normalized by $\Sc_n$. Thus the semi-direct product of $\Sc_n$ and the group of sign changes yields a reflection group $W= (\ZZ / 2 \ZZ)^n \rtimes  \Sc_n$, called the real reflection group of type $B_n$, denoted by $W(B_n).$\\
 
 A root system $ R$ associated to a reflection group $W$ is a finite set  of nonzero vectors in $\fh$ satisfying the five following conditions: 
 \begin{enumerate}
\item $R$ spans $\fh$ as a vector space;
\item $R \cap \bR \alpha =\{ \alpha, -\alpha \}$ for all $\alpha \in R;$
\item $s_\alpha( R)= R$ for all $\alpha \in R$;
\item for any $\alpha, \beta \in R$, $\displaystyle \frac{(\alpha, \beta)}{(\beta, \beta)}\in \ZZ $;
\item $W$ is generated by all reflections $s_\alpha,  \alpha \in R.$
 \end{enumerate}
 Let $\eps_1, \ldots,\eps_n$ the standard basis of $\fh=\bR^n$. The set corresponding set of roots $R$ associated with $W(B_n)$ consists of $2n$  roots $\pm \eps_i$ and and $2n(n-1)$   roots $\pm \eps_i \pm \eps\ (i <j)$ \cite{Hu}.

\subsection{ Quantum Olshanetsky-Perelomov Hamiltonian }

Let $W $ be a real reflection group, $R$ the system of roots associated with $W$  and  
$ S= \{ s_\alpha | \alpha \in R \}$ the set of reflections. 
Clearly, $W$ acts on $S$ by conjugate. Let $c : S \to \bC$ be a conjugate invariant function.
The quantum Olshanetsky-Perelomov Hamiltonian attached to $W$ is the second order differential operator
$$ H := \Delta_{\fh} - \sum_{s \in S}  \frac{c_s(c_s+1)(\alpha_s, \alpha_s)}{\alpha_{s}^2},$$
where $\Delta_{\fh}$ is the Laplace operator on $\fh$.
It turns out that the system defined by the Olshanetsky-Perelomov operator $H$ is completely integrable. Namely, we have the following theorem.
\begin{theorem} \cite{OP}
There exist differential operators $L_j$ on $\fh$ with  rational coefficient  such that $L_j$ are homogeneous (of degree $-d_j$),\ $L_1 =H,$ and $[L_j,L_k]=0,\ \forall j, k.$
\end{theorem}This Theorem is obviously a generalization of \cite[Theorem 2.1]{DOS}. \\
Let $y \in \fh$, The {\it Dunk-Opdam operator} $D_y$ on $\bC(\fh)$ is defined by the formula  \cite{Dkl-Opd}

$$ D_y:= \partial_y - \sum_{s \in S} c_s \frac{ ( \alpha_s, y)  }{\alpha_s} (1-s)$$ 
where $\partial_y$ denotes the partial derivatives in the direction of $y$. 
Clearly, $D_y \in  \D (\fh_{reg})   \rtimes \bC[W]$, where $\fh_{reg}$ is the set of regular points of $\fh$ and $\Dc(\fh_{reg})$ denotes the algebra of differential operators on $\fh_{reg}$.

 For any element $ B \in \D (\fh_{reg})   \rtimes \bC[W],$ define $m(B)$ to be the differential operator $\bC(\fh)^W \to \bC(\fh),$ attached to $B.$ That is, if $B= \sum_{g \in W} B_g g,\ B_g \in \D (\fh_{reg}),$  then $m(B)=\sum_{g \in W} B_g.$ It is clear that if $B$ is  $W$-invariant, then for all  $A \in \bC[W] \ltimes \D (\fh_{reg}),$
 $$ m(AB)= m(A) m(B).$$
Let $S\hf$ be  symmetric algebra of $\hf$ and $(S\fh)^W= \{ x \in S\fh |\ wx=x, \forall w \in W \}$. Let us recall   by the Chevalley-Shepard-Todd theorem, that  the algebra $(S\hf)^W$ is free \cite{Che}. Let $P_1, \ldots, P_r$ be homogeneous generators of $(S\fh)^W$.
\begin{proposition}\cite{Du1}
Let $\{ y_1, \ldots, y_r \}$ be an orthonormal basis of $\fh$, and  $P_1, \ldots, P_r$ be homogeneous generators of $(S\fh)^W$. Then we have
$$ L_j= m(P_j(D_{y_1}, \ldots, D_{y_r})).$$ 
\end{proposition}
It is  now clear that the  differential operators  $H$ and the $L_j$ are invariant under the action of the reflection group $W$ so they belong to a localization of the ring of invariant differential operators under the reflection group $W$. \\  The action of  the Olshanetsky-Perelomov operator $H$ on polynomials could be understood by studying the polynomials representation of the ring of  differential operators  invariant under the reflection group $W$ localized at  $\Delta^2,$ where $\displaystyle \Delta= \prod_{s \in S} \alpha_s(x)$ is the polynomial discriminant.\\

\subsection{Higher Specht Polynomials for Reflections group $G(r,p)$} In this subsection we recall  some general facts about the representation of wreath product  $G(r,n)$. Let $\Sc_n$ be the group of permutations of the set of variables  $\{ x_1, \ldots, x_n\}$ and $\ZZ /r \ZZ$ be the cyclic group of order $r$ which acts on the  $x_i$ by a primitive $r$th root of unity.\\
The  wreath product $G(r,n)$ is  the semi-direct product of $(\ZZ/r \ZZ)^n$ with $\Sc_n$, written as $(\ZZ / r\ZZ)^n \rtimes \Sc_n$, where $(\ZZ / r\ZZ)^n $ is the direct product of $n$ copies of $\ZZ / r\ZZ $. Let $\xi$ be a primitive $r$-th root of  1.
$(\ZZ / r\ZZ)^n \rtimes \Sc_n= \{ (\xi^{i_1}, \ldots, \xi^{i_n}; \sigma) |\  i_k\in \bN,\  \sigma \in \Sc_n \}$, whose product is given by 
$$  (\xi^{i_1}, \ldots, \xi^{i_n}; \sigma) (\xi^{j_1}, \ldots, \xi^{j_n}; \pi)=  (\xi^{i_1 + j_{\sigma^{-1} (1)}}, \ldots, \xi^{i_n + j_{\sigma^{-1} (n)}},; \sigma \pi).$$
Let $\Oc_X= \bC[x_1, \ldots,x_n]$ be the ring of polynomials in $n$ indeterminates  on which the group $G(r,n)$ acts as follows:

$$ (\xi^{i_1}, \ldots, \xi^{i_n}; \sigma) f= f(\xi^{i_{\sigma(1)}} x_{\sigma(1)}, \ldots, \xi^{i_{\sigma(n)}}x_{\sigma(n)}; \sigma),$$
where $f \in \Oc_X$ and  $(\xi^{i_1}, \ldots, \xi^{i_n}; \sigma) \in G(r,n)$. It is known that the fundamental invariants under this action are given by the elementary symmetric functions $e_j(x_1^r, \ldots, x_n^r),\ 1 \leq j \leq n.$ Let $J_{+}$ be the ideal of $\Oc_X$ generated by these fundamental invariants and $\Lambda=\Oc_X/ J_{+}$ be the quotient ring. It is also known that the $G(r,n)$-module $\Lambda$ is isomorphic to the group ring $\bC[G(r,n)]$, namely the left regular representation. A description of all irreducible components of $\Lambda$ is known in \cite{Ariki}, in terms of what is called  "higher Specht polynomials". The irreducible representation of $G(r,n)$ are parametrized by the $r$-tuple of Young diagrams $(\lambda^1,\ldots, \lambda^r)$ with $|\lambda^1| +\cdots + |\lambda^r|=n.$ Let $\mathcal{P}_{r,n}$ be the set of $r$-tuples of Young diagrams $\bla=(\lambda^1, \ldots, \lambda^r)$ with $|\lambda^1| +\cdots + |\lambda^r|=n$. By filling each cell with a positive integer in such a way that every  $j\ ( 1 \leq j \leq n)$ occurs once, we obtain an $r$-tableau $T=(T^1,\ldots,T^r)$ of shape ${\bf \bla}= (\lambda^1, \ldots, \lambda^r)$. When the number $k$ occurs  in the component $T^i$, we write $k \in T^i$. The set of $r$-tableaux of shape $\bla$ is denoted by $\Tab(\bla).$ An $r$-tableau $T=(T^1,\ldots,T^r)$ is said to be standard if the numbers are increasing on each column and each row of $T^\nu\ (1 \leq \nu \leq  r)$. The set of $r$-standard tableaux of shape $\bla$ is denoted by $\STab(\bla).$

Let $S =(S^1, \ldots, S^r) \in \STab (\bla)$. We associate a word $w(S)$ in the following way. First we read each column of the component $S^1$ from the bottom to the top starting from the left. We continue this procedure for the tableau $S^2$ and so on. For word $w(S)$ we define index $i(w(S))$ inductively as follows. The number 1 in the word $w(S)$ has the index $i(1)=0.$ If the number $k$ has index $i(k)=p$ and the number has number $k+1$ is sitting on the left (resp. right) of $k,$ then $k+1$ has index $p+1$ (resp. $p$). Finally, assigning the indices to the corresponding cells, we get a shape ${\bf \bla}= (\lambda^1, \ldots, \lambda^r)$, each cell filled with a nonnegative integer, which is denoted by $i(S) =(i(S)^1, \ldots, i(S)^r).$

Let $T=(T^1,\ldots, T^r)$ be an $r$-tableau of shape $\bla$. For each component $T^\nu\ (1\leq \nu \leq r),$ the Young symmetrizer $\eb_{T^\nu}$ of $T^\nu$ is defined by
$$ \eb_{T^\nu}= \frac{1}{ \alpha_{T^\nu}} \sum_{\sigma \in R(T^\nu)\ \tau \in C(T^\nu)} \sgn (\tau) \tau \sigma,$$ where $\alpha_{T^\nu}$ is the product of the hook lengths for the shape $\lambda^\nu$,  $R(T^\nu)$ and $C(T^\nu)$ are the \textit{row-stabilizer} and \textit{colomn-stabilizer} of $T^\nu$ respectively.\\
We may regard a tableau $T$ on a Young diagram $\lambda$ as a map
$$ T: \{ \mbox{cells\ of}\  \lambda \}  \to \ZZ_{\geq 0},$$
which assigns to a cell $\xi$ of $\lambda$ the number $T(\xi)$ written in the cell $\xi$ in $T$.\\
For $S\in \STab(\bla)$ and $T \in\Tab(\bla)$, Ariki, Terasoma and Yamada in \cite{Ariki} defined the higher Specht polynomial for $G(r,n)$ by
$$ F_{T}^{S}= \prod_{\nu =1}^r \bigg( \eb_{T^\nu} (x_{T^\nu}^{ri(S)^{\nu}}) \prod_{k\in T^\nu} x_k^\nu \bigg),$$ where $$\displaystyle x_{T^\nu}^{ri(S)^{\nu}}= \prod_{\xi \in \lambda^\nu} x_{T^\nu(\xi)}^{ri(S)^\nu(\xi)}.$$

The following is the fundamental result in \cite{Ariki} on the higher Specht polynomials for $G(r,n).$
\begin{theorem}
\begin{enumerate}
\item The space $\displaystyle V_S(\lambda) =\sum_{T\in \Tab(\lambda)} \bC F_T^S$ affords an irreducible representation of the reflection group $G(r,n)$.
\item The set $\{ F_T^S \ | \ T \in \STab(\bla) \}$  gives a basis over $\bC$ for $V_S(\bla).$ 
\item For $S_1 \in \STab(\bla)$ and $S_2 \in \STab(\mu)$, the representation $V_{S_1}(\bla)$ and $V_{S_2}(\mu)$ are isomorphic if and only if $S_1$ and $S_2$ has the same shape, i.e. $\lambda=\mu.$

\item We have the irreducible decomposition 
$$ \bC[G(r,n)]= \bigoplus_{\bla \in \mathcal{P}_{r,n}} \bigoplus_{S\in \STab(\bla)} V_S(\lambda)$$
as representation of $G(r,n)$.
\end{enumerate}
\end{theorem}

\begin{theorem}

The higher Specht polynomials in $\C=\{ F_T^S; S, T \in \STab(\lambda), \lambda \vdash n \}$
form a basis of the $\bC[x_1,...,x_n]^{G(r,n)}$-module $\bC[x_1,...,x_n]$
\end{theorem}

\section{Decomposition Theorem} 
In this section we establish a decomposition theorem of the polynomial ring in $n$ indeterminates as a module over the ring of invariant differential operators.

We are interested in studying the action of the invariant differential operators under the real  reflection group  $W=W(B_n)$ of type $B_n$. We know  that $W(B_n)= (\ZZ / 2\ZZ)^n \rtimes \Sc_n$. Let $\eps_1, \ldots,\eps_n$ the standard basis for $\fh=\bR^n$. The associated set of roots $R$ consists of $2n$ short roots $\pm \eps_i$ and and $2n(n-1)$  long roots $\pm \eps_i \pm \eps\ (i <j)$. Then the polynomial discriminant $\displaystyle \Delta= \prod_{s \in S} \alpha_s(x) = 2^n n! x_1\cdots x_n \prod_{1\leq i<j \leq n} (x_j^2-x_i^2).$

\subsection{Actions description}
As we want to study the polynomial representation of a ring of invariant differential operators localized at $\Delta^2,$ it is convenient to precisely describe the action of that ring  on the polynomials ring.\\

Let $\Dc_X=\bC \langle x_1, \ldots, x_n, \frac{\partial}{\partial x_1}, \ldots, \frac{\partial}{\partial x_n} \rangle$ be the ring of differential operators associated with the polynomial ring $\Oc_X=\bC [x_1,\ldots,x_n]$, and $\Oc_Y=\bC[x_1, \ldots,x_n]^{W}= \bC[y_1,\ldots,y_n]$ be the ring of invariant polynomials under the real reflection group $W$ where  $$y_j =\displaystyle \sum_{i=1}^n {x_i^{2j}} \  \mbox{for} \  j=1,\ldots, n,$$

We denote by $\Dc_Y= \bC \langle y_1, \ldots, y_n, \frac{\partial}{\partial y_1}, \ldots, \frac{\partial}{\partial y_n} \rangle $ the ring of differential operators associated with  $\Oc_Y= \bC [y_1,\ldots,y_n]$.  By \cite{LS1}, $\Dc_Y$ is the ring of invariant  differential operators under the action of the real reflection group $W$. It is  not  clearly that $\Oc_X$ is a $\Dc_Y$-module, we need to describe the action of $\Dc_Y$ on $\Oc_X$.  By localization, $\Oc_X$ is turned into a $\Dc_Y$-module, as the following lemma states. \\

\textbf{Notations} We adopt the following notations  $$\displaystyle \OO:=\bC [x_1,\ldots,x_n, \Delta^{-1}],\  \tilde{\Oc_Y} := \bC [y_1,\ldots,y_n, \Delta^{-2}],$$ $$  \DD:= \bC \langle y_1, \ldots, y_n, \frac{\partial}{\partial y_1}, \ldots, \frac{\partial}{\partial y_n}, \Delta^{-2} \rangle.$$

\begin{lemma}

$\tilde{\Oc_X}$ is a $\tilde{\Dc_Y}$-module. 
\end{lemma}
\begin{proof} 
Let us make clear the  action of $\tilde{\Dc_Y}$ on $\tilde{\Oc_X}.$\\
We have    $y_j =\displaystyle \sum_{i=1}^n {x_i^{2j}}, j=1,\ldots, n,$  hence  $\displaystyle  \frac{ \partial}{\partial x_i} = \sum_{j=1}^n 2j x_i^{2j-1} \frac{ \partial }{ \partial y_j} , \ i= 1,\ldots, n.$ Let $A=(x_i^{2j-1})_{1\leq i,j \leq n}$ so that $\det(A)= \Delta$. We get the following equation
$$ \left (\begin{array}{ccc} \frac{\partial}{\partial x_n}\\  \vdots\\ \frac{\partial}{\partial x_n} \end{array} \right)= A  \left (\begin{array}{ccc} \frac{\partial}{\partial y_1} \\ \vdots \\  \frac{\partial }{\partial y_n} \end{array} \right).$$
Since $\Delta \neq 0$, it follows that 
$$ 
\left (\begin{array}{ccc} \frac{\partial}{\partial y_1} \\ \vdots \\  \frac{\partial }{\partial y_n} \end{array} \right)= A^{-1} \left (\begin{array}{ccc} \frac{\partial}{\partial x_n}\\  \vdots\\ \frac{\partial}{\partial x_n} \end{array} \right)
$$
and it is now clear that $\tilde{\Oc_X}$ is a $\tilde{\Dc_Y}$-module. 
\end{proof}

{\it  Is $\OO$ a $\DD$-  semisimple module ? If yes what are the simple components of $\OO$ as $\DD$-module  and their multiplicities?}\\

\subsection{Simple components and their multiplicities}

In this section, we state our main result. We use the representation theory of the real reflection group $W$ to yield results on modules over the ring of differential operators. It is well-known that 
$$ \Oc_X= \bC[W] \otimes \Oc_Y \  \mbox{as } \  \Oc_Y\mbox{-modules}.$$ 
Let us consider the multiplicative closed set $S=\{ \Delta^{k} \}_{k \in \bN} \subset \Oc_X.$ It follows that:
$$ S^{-1} \Oc_X= \bC[W] \otimes S^{-1}\Oc_Y \  \mbox{as } \  S^{-1}\Oc_Y\mbox{-modules}.$$ where $S^{-1} \Oc_X$ and $S^{-1}\Oc_Y$ are the localizations of  $\Oc_X$ and $\Oc_Y$ at $S$ respectively. But $S^{-1} \Oc_X= \OO$ and $S^{-1} \Oc_Y= \tilde{\Oc_Y}$ , whereby we get
$$ \OO=  \bC[W] \otimes \tilde{\Oc_Y} \  \mbox{as} \  \bC[W]\mbox{-modules}.$$

\begin{lemma}
There exists an injective map
$$ \bC[W] \hookrightarrow  \Homo_{\bC}( \OO, \OO ).$$
\end{lemma}
\begin{proof}
The $\bC[W]$-module  $\bC[W]$ acts on itself by multiplication, and this multiplication yields an injective map $\bC[W] \hookrightarrow \Homo_{\bC} ( \bC[W] , \bC[W]).$ Since $\tilde{\Oc_Y} $ is invariant under this action of $\bC[W]$, we get the expected injective map. 
\end{proof}

\begin{proposition}
There exists an injective map 
$$ \bC[W] \hookrightarrow  \Homo_{\tilde{\Dc_Y}}( \OO , \OO ).$$
\end{proposition}
\begin{proof}
Since  $\tilde{\Dc_Y}= \bC \langle y_1, \ldots, y_n, \partial y_1, \ldots, \partial y_n, \Delta^{-2}  \rangle $, we only need to show that every element of $\bC[W]$ commutes with $y_1,\ldots,y_n, \partial y_1, \ldots, \partial y_n.$ 
\begin{enumerate}
\item[$\bullet$] It is clear that every element of $\bC[W]$ commutes with $y_i,\ i=1, \ldots,n$.
\item[$\bullet$] Let us show that every element of $\bC[W]$ commutse with $\partial y_i,\ i=1, \ldots,n$.
Let ${\bf D}$ be a derivation on the field $ K=\bC(y_1, \ldots, y_n)$  of fractions of $\Oc_Y$, then $(K, {\bf D})$ is a differential field. Let $L=\bC(x_1, \ldots, x_n)$ be the field of fractions of $\OO$. We have that  $K=L^W$ is the fixed field and $L$ is a Galois extension of $K$, with Galois group $W$.  Then by \cite[Th\'eor\`eme 6.2.6]{Cham}  there exists a unique derivation on $L$ which extends ${\bf D},$ then  $(L , {\bf D})$ is also a differential ring. In this  way, $\sigma^{-1}  {\bf D} \sigma = {\bf D}$ for every $\sigma \in W.$ Therefore $\sigma {\bf D}= {\bf D}\sigma$ and $ \sigma $ commute with ${\bf D}.$   
\end{enumerate}
\end{proof}

\begin{corollary}
$$\bC[ W] \cong \Homo_{\tilde{\Dc_Y}} ( \OO , \OO )$$ 
\end{corollary}
\begin{proof}
see \cite[Corollary 26 ]{IBK}
\end{proof}
Before we state our  main result, let us recall some facts.\\
By Maschke's Theorem \cite[Chap XVIII]{Lan}, we know that $\bC[W]$ is a semi-simple ring, and  

$$\bC[W] = \displaystyle \bigoplus_{ \lambda \in \Pc_{2,n}} R_{\lambda}, $$

 where  $\mathcal{P}_{2,n}$ be the set of $2$-tuples of Young diagrams $\bla=(\lambda^1,  \lambda^2)$ with $|\lambda^1| + |\lambda^2|=n$
  and  $R_{\lambda}$ are simple rings. In fact $ \displaystyle R_\bla = \bigoplus_{ S \in \STab(\bla) } V_S (\bla) $ ( see  Theorem 2.3).  
  We have the following corresponding decomposition of the identity element of $\bC[W]$:
$$ 1=\sum_{ \lambda \in \Pc_{2,n}}  r_{\lambda},$$
 where $r_{\lambda}$ is the identity element of $R_{\lambda}$, with  $r_{\lambda} ^2=1$ and $r_{\lambda} r_{\mu} =0$ if $\lambda \neq \mu.$  $\{ r_{\lambda} \}_{ \lambda \in \Pc_{2,n}}$ is the set of primitive central idempotents of $\bC[W]$.

Let $n \in \bN^*, \bla \in  \Pc_{2,n}$, we set  $ \displaystyle \Tab(n)=\displaystyle \cup_{\lambda \in \Pc_{2,n}} \Tab(\lambda)$ and $\displaystyle \STab(n)=\cup_{\lambda \in \Pc_{2,n}} \STab(\lambda)$.

\begin{theorem}

For every  primitive idempotent $e \in \bC[W]$.
\begin{enumerate}
\item $ e \OO$ is a nontrivial  $\tilde{\Dc}_Y$-submodule of $\OO,$
\item The $\DD$-module $e \OO$ is simple,
\item There exist  $\lambda  \in \Pc_{2,n}$ and a higher Specht polynomial $F_T^S$  (with $S,T \in \STab(\lambda))$ such that $e  \OO= \DD F_T^S.$
\end{enumerate}
\end{theorem}
\begin{proof}\
\begin{enumerate}
\item Let $e  \in \bC[W]$ be a  primitive  idempotent, we know that $\bC[W]e $ is a $W$- irreducible representation.  Theorem 2.3 states that  there is $\lambda \in \Pc_{2,n}$ and $S\in \STab(\lambda)$ such  that $\bC[W] e \cong V_S(\lambda)$, and $V_S(\lambda) \subset \OO.$ By  \cite[ Chap III, \S 4, Theorem 3.9 ]{Boerner}, we have $ e \bC[W]e \cong \bC e  \neq \{0\}$. $ \{0\} \neq e V_S(\lambda) \subset e \OO.$  Since $e $ commute with every element of $\DD$ and $\OO$ is a $\DD$-module, it follows that $e  \OO$ is a nontrivial $\DD$-module.\\
In fact $V_S(\lambda)$ is a cyclic $\bC[W]$-module, i.e.,  there exist  $T, S \in \STab(\lambda)$ and a higher Specht polynomial $F_T^S$ such that  $ V_S(\lambda) = \bC[W]  F_T^S$, so that $ \bC[W] e \cong \bC[W] F_T^S$. Then it follows that  that $e  F_T^S$ is a scalar multiple of $F_T^S.$
 
\item Assume that $1 =\sum_{i=1}^s e_i$ where the $\{ e_i \}_{ 1 \leq i \leq s}$ is the set of primitive idempotents of $\bC[W],$ then $\OO= \sum_{i=1}^s e_i \OO.$ Let $m \in e_i \OO \cap e_j\OO$  with $ i \neq j$ so that $m=e_i m_i$ and $ m= e_j m,$ but $ e_i e_j=0$ then $ e_im=e_ie_jm=0$ hence $m=0.$ Therefore 
$ \OO = \oplus_{i=1}^s e_i \OO$  {and}  we get:   
$$ \Homo_{\DD} (\OO, \OO) \cong  \bigoplus_{i,j =1}^s   \Homo_{\DD} (e_i\OO, e_j\OO), $$
by Corollary 3.4   we know that
$ \displaystyle \bC[W] \cong \bigoplus_{i,j =1}^s   \Homo_{\DD} (e_i\OO, e_j\OO).$
For every $\lambda \in \Pc_{2,n}$, we pick a unique  irreducible representation $V_S(\lambda)$ for a certain $S\in \STab(\lambda)$ which we denote by $V(\lambda):= V_S(\lambda)$.
We also have, by \cite[Proposition 3.29]{Fulton-Harris}, that  
$ \displaystyle \bC[W] \cong \bigoplus_{\lambda \in \Pc_{2,n} } \mbox{End}_{\bC} (V({\lambda})).$  
But by the Wedderburn's decomposition Theorem \cite[Chap II,\S 4, Theorem 4.2]{Boerner} we also know that 
$$ \displaystyle \bC[W] = \bigoplus_{\lambda \in \Pc_{2,n}} r_{\lambda} \bC[W]\ \mbox{and}\ r_{\lambda} \bC[W]  \cong \mbox{Mat}_{f^{\lambda}}(\bC) \cong \mbox{End}_{\bC} ( \bC^{f^{\lambda}}) $$ {where}  $f^{\lambda} = \dim_{\bC} (V(\bla)).$ We recall that each primitive idempotent $e_i$   is associated  with standard tableau $T_i,$ we may denote  $e_i=e_{T_i}$.Let us show that  $$ \displaystyle \bC[W] \cong  \bigoplus_{\lambda \vdash n} \bigg( \bigoplus_{T_i, T_j \in \STab(\lambda)}   \Homo_{\DD} (e_i\OO, e_j\OO) \bigg)\  \mbox{where}\ e_i=e_{T_i}.$$ \\
 Let $x$ be an element of $ \bC[W]$ and $r_{\lambda}$ the primitive central idempotent associated with $\lambda \in \Pc_{2,n}.$ Then $x$ induces an $\DD$-homorphism $\OO \to \OO, m \mapsto x\cdot m;$ the multiplication by $x$. Since $r_{\lambda}$ is in the centre of $\bC[W] $,  $ x \cdot  (r_{\lambda} \OO )=(x \cdot r_{\lambda} ) \OO \subset r_{\lambda} \OO,$ which means $x \in \displaystyle \oplus_{\lambda \in \Pc_{2, n}} \Homo_{\DD} ( r_{\lambda} \OO, r_{\lambda} \OO).$  It follows that
 $$ \displaystyle \Homo_{\DD} (  \OO, \OO)  \cong  \bigoplus_{\lambda \in \Pc_{2,n}}    \Homo_{\DD} (r_\la \OO, r_\la\OO) .$$
  Then $\Homo_{\DD} ( e_i \OO , e_j \OO) =\{0 \}$ if $T_i \in \STab(\lambda_i), T_j \in \STab(\lambda_j)$ and  $\lambda_i \neq \lambda_j.$ 
  We  get that   $$ \displaystyle \Homo_{\DD} (  \OO, \OO)  \cong  \bigoplus_{\lambda \in \Pc_{2,n}} \bigg( \bigoplus_{T_i, T_j \in \STab(\lambda)}   \Homo_{\DD} (e_i\OO, e_j\OO) \bigg).$$ 
 The number of direct factors in the sum \\$\displaystyle  \bigoplus_{T_i, T_j \in \STab(\lambda)}   \Homo_{\DD} (e_i\OO, e_j\OO) $ is $(f^{\lambda})^2.$\\
Let us show that $  \Homo_{\DD} (e_i\OO, e_j\OO) \cong \bC$ if $T_i, T_j \in \Tab(\lambda).$  Consider the following commutative diagram:

\begin{center}
\begin{tikzpicture}
  \matrix (m) [matrix of math nodes,row sep=3em,column sep=4em,minimum width=2em]
  {
    \bC[W] & \displaystyle \Homo_{\DD} (  \OO, \OO)  \\
     r_{\lambda} \bC[W]& \displaystyle \Homo_{\DD} (  r_{\lambda}\OO, r_{\lambda}\OO) \\};
  \path[-stealth]
    (m-1-1) edge node [left] {$\alpha_{\lambda}$} (m-2-1)
            edge  node [below] {$\phi$} (m-1-2)
    (m-2-1.east|-m-2-2) edge node [below] 
    {$\psi$}
            node [above]{ }(m-2-2)
    (m-1-2) edge node [right] {$\beta_{\lambda}$} 
    (m-2-2); edge [dashed,-] (m-2-1);
\end{tikzpicture}
\end{center}

where
 $\displaystyle \beta_{\lambda}: \bigoplus_{\mu \in \Pc_{2,n}} \Homo_{\DD} ( r_{\mu} \OO, r_{\mu} \OO) \to 
 \Homo_{\DD} ( r_{\lambda} \OO, r_{\lambda}\OO)$ and \\ $ \displaystyle \alpha_{\lambda}: \bigoplus_{\mu \in \Pc_{2, n}} r_{\mu} \bC[W] \to r_{\lambda} \bC[W]$ are canonical projections et $\phi$ is the isomorphism in Corollary3.4. It follows that $\psi$ is an isomorphism hence $ r_{\lambda} \bC[W] \cong   \displaystyle \Homo_{\DD} (  r_{\lambda}\OO, r_{\lambda}\OO).$ \\
 Now  we  identify $ r_{\lambda} \bC[W]$ with either the set  $\mbox{Mat}_{f^{\lambda}}(\bC)$ of square matrices of order $f^{\lambda}$ with coefficients in $\bC$ either with  $\displaystyle \mbox{End}_{\bC}(\bC^{f^{\lambda}}).$
 Let $E_{ij}$ be the square matrix of order $f^{\lambda}$  with 1 at the position $(i,j)$ and 0 elsewhere and $E_i=E_{i,i}, $ then we identify the primitive  idempotent $e_i \in  r_{\lambda} \bC[W]$ with $E_i$ in $\mbox{Mat}_{f^{\lambda}} (\bC).$ Let $B=(a_{ij}) \in \mbox{Mat}_{f^{\lambda}}(\bC)$ we get $B=\sum_{i,j} a_{i,j} E_{i,j}= \sum_{i,j} E_i B E_{j}$, in fact $E_iBE_j$ is  the matrix with $a_{i,j}$ in the position $(i,j)$ and 0 elsewhere, if $R= \mbox{Mat}_{f^{\lambda}}(\bC)$  we get that  $ E_i R E_j \cong \bC.$ 
 
  This isomorphism $\psi$ implies that $$ \displaystyle  \bigoplus_{T_i,T_j \in \STab(\lambda)} E_i R E_j \cong  \bigoplus_{ T_i,T_j \in \STab(\lambda)} \Homo_{\DD} (e_i\OO, e_j \OO);$$ the restriction of $ \psi$ to $E_iRE_j$ yields a map\\ $E_i R E_j \to \Homo_{\DD}(e_i \OO, e_j \OO)$ and this map is surjective, moreover we have $E_i R E_j \cong \Homo_{\DD}(e_i \OO, e_j \OO)$. \\Therefore $ \Homo_{\DD}(e_i \OO, e_i \OO) \cong \bC.$ Let us assume that $e_i \OO$ is not  simple $\DD$-module, then $e_i\OO$ may be written as $ e_i \OO= \oplus_{j \in J} N_j$ where the $N_j$ are simple $\DD$-modules and $|J| > 1$. It follows that $ \dim_{\bC} (  \Homo_{\DD}(e_i \OO, e_i \OO))  \geq |J|$  but \\ $ \Homo_{\DD}(e_i \OO, e_i\OO) \cong \bC$ so we obtain that $J=1$, which necessary implies that $e_i \OO$ is a simple $\DD$-module. 
 \item By the the proof (i) there exists a higher Specht polynomial $F_T^S  \in e_i \OO$, with $ S, T_i\in \STab(\bla) \lambda \vdash n$ such that $ e_i \OO = \DD F_{T_i}^S.$
\end{enumerate}
\end{proof}

\begin{corollary}

With the above notations, $e_i\OO \cong_{\DD} e_j \OO$ if only if $T_i$ and $T_j$  have the same shape i.e. if there is a partition $\lambda \in \Pc_{2,n} $ such that $ T_i, T_j \in \STab(\lambda ).$
\end{corollary}

\begin{proof}
The  $\DD$-modules $e_i \OO$ are simple and $\Homo_{\DD}(e_i \OO, e_j \OO) \cong \bC$ whenever there exists a partition $\lambda \in \Pc_{2,n}$ such that $ T_i, T_j \in \STab(\lambda).$ Since $  \Homo_{\DD}(e_i \OO, e_j \OO)\neq \{0\},$ we conclude by using the Schur lemma.

\end{proof}

\begin{proposition}
  Let  $\lambda \in \Pc_{2,n}$, $ T \in \STab(\lambda) $, and let $e$ be the  primitive idempotent associated  with $T,$  denote by $F_T:=F_T^S$  the corresponding higher Specht polynomial (for some $S\in \STab(\bla)$), in Theorem 3.5 (iii),  such that $e\OO=\DD F_T^S$  then  we have: 
  \begin{enumerate}
\item
 \begin{equation}
  \displaystyle \OO = \bigoplus_{{T\in \STab(n)}}  \DD F_T= \bigoplus_{\bla \in \Pc_{2,n}} \bigg( \bigoplus_{{T\in \STab(\bla)}}  \DD F_T \bigg);
  \end{equation}
\item for each  $\lambda \in \Pc_{2,n}$  fix  a 2-tableau $ T^* \in \STab(\lambda) $, then 
 \begin{equation}
  \displaystyle \OO =  \bigoplus_{\bla \in \Pc_{2,n}} f^\bla \DD F_{T^*}
  \end{equation}
  where $f^\bla =\dim_{\bC} (V(\lambda))$
\end{enumerate}

\end{proposition}
\begin{proof}
We have by the proof of  Theorem 3.5  that
$$ \OO = \bigoplus_{T_i \in \STab(n)} e_{i} \OO $$ and the $e_i \OO$ are simple $\DD$-modules. Since to each primitive idempotent $e_i$ corresponds a $2$-diagram  $\lambda_i \in \Pc_{2,n} $ and a tableau $T_i \in \STab(\lambda_i)$ such that $e_i\OO= \DD F_{T_i}$ then $\OO= \bigoplus_{T \in \STab(n)} \DD F_T.$ By Corollary 3.6, $ \DD F_{T_j} \cong  \DD F_{T_j}$ if $T_i,T_j \in \STab(\lambda)$ and so we have $f^\bla$ isomorphic copies of $ \DD F_{T^*}$ in the direct sum (3.1).

\end{proof}

We get in Proposition 3.7 a decomposition of the polynomial ring as a $\DD$-module into irreducible $\DD$ modules generated by the higher Specht polynomials.

\section*{Acknowledgments}
The final version work has been done while the first author was visiting IMSP at Benin,  under the  {\it Staff Exchange} program of the German  Office of Academic Exchange  (DAAD).
He warmly thanks   the DAAD for the financial support.


\end{document}